\documentclass[12pt]{article}
\usepackage[utf8]{inputenc}
\usepackage{amssymb,amsmath,amsfonts,eucal,mathrsfs,amsthm} 
\usepackage[colorinlistoftodos]{todonotes}
\usepackage{enumitem}
\newtheorem{theorem}{Theorem}
\newtheorem{proposition}[theorem]{Proposition}
\newtheorem{lemma}[theorem]{Lemma}

\newtheorem{corollary}[theorem]{Corollary}

\newtheorem{remark}[theorem]{Remark}

\newtheorem{example}[theorem]{Example}

\theoremstyle{definition}
\newcommand{\R}{\mathbb{R}}

\newcommand{\Sf}{\mathbb{S}}

\newcommand{\Hy}{\mathbb{H}}
\newcommand{\spa}{\mbox{span}}
\newcommand{\hess}{\mbox{Hess\,}}
\newcommand{\Ric}{\mbox{Ric}}

\newcommand{\nab}{\tilde\nabla}

\newcommand{\U}{\cal{U}}
\newcommand{\p}{\partial}
\def\<{{\langle}}
\def\>{{\rangle}}

\def\n{\nabla}
\def\d{\partial}
\def\a{\alpha}

\def\be{\begin{equation} }
\def\ee{\end{equation} }

\newcommand\blfootnote[1]{
\begingroup
\renewcommand\thefootnote{}\footnote{#1}
\addtocounter{footnote}{-1}
\endgroup
}
\begin{document}

\title{A class of Einstein submanifolds of\\ Euclidean space}
\author{M. Dajczer, C.-R. Onti and Th. Vlachos}
\date{}
\maketitle
\blfootnote{\textup{2020} \textit{Mathematics Subject Classification}:
53C25,  53C40, 53C42.}

\blfootnote{\thanks{This research is in part a result of the activity 
developed within the framework  of the Programme in Support of Excellence 
Groups of the Regi\'on de Murcia, Spain, by Fundaci\'on S\'eneca, Science 
and Technology Agency of the Regi\'on de Murcia.
Marcos Dajczer was partially supported by MICINN/FEDER project
PGC2018-097046-B-I00 and Fundaci\'on S\'eneca project 19901/GERM/15, Spain.}}

\begin{abstract} In this paper we give local and global parametric
classifications of a class of Einstein submanifolds of  Euclidean 
space. The highlight is for submanifolds of codimension two since 
in this case the assumptions are only of intrinsic nature.
\end{abstract}

Let $f\colon M^n\to\R^{n+p}$ denote an isometric immersion of an 
$n$-dimensional Riemannian manifold into Euclidean space with 
codimension $p$. The goal of this paper is to classify parametrically 
certain classes of these submanifolds for which $M^n$ is an 
Einstein manifold of non-constant sectional curvature. Recall 
that a Riemannian manifold $M^n$ is said to be \emph{Einstein} if 
its Ricci tensor is proportional to the metric, that is, if 
$$
\Ric_M(X,Y)=\rho\,\<X,Y\>
$$
for any vector fields  $X,Y\in\mathfrak{X}(M)$ and some constant 
$\rho\in\R$. Hence $\rho$ is the (not normalized) constant Ricci 
curvature of $M^n$.

We focus on manifolds of dimension $n\geq 4$ since $3$-dimensional 
Einstein manifolds have constant sectional curvature. Besides, 
we are interested in codimension $p\geq 2$ since, due to an observation 
by Cartan communicated by Thomas \cite{Th} in 1937 and the work of 
Fialkow \cite{Fi} from 1938, we know 
that an Einstein hypersurface $f\colon M^n\to\R^{n+1}$, $n\geq 3$, 
is either flat or an open subset of a round sphere; see \cite{DT}
for a proof. In particular, if $M^n$ is a complete manifold, then 
the submanifold is either a cylinder over a complete plane curve 
or a round sphere.

At the present time, the knowledge on the subject of Euclidean 
Einstein submanifolds, except those with constant sectional curvature, 
is quite limited. This is because the assumption that a 
submanifold is Einstein is much weaker than having constant sectional 
curvature. In fact, as far as we know, until now  the only classification 
result available under purely intrinsic assumptions is the
aforementioned for hypersurfaces.
Hence, in order to obtain a meaningful classification one has to require
some strong additional hypothesis. For instance, a classification was
obtained by Onti \cite{Ont} under the extrinsic assumptions that the 
submanifold has parallel mean curvature vector and flat normal bundle.

It has been shown by Dajczer, Onti and Vlachos \cite{DOV} that 
any Einstein Euclidean submanifold with flat normal bundle is 
locally holonomic. This means that the manifold carries a system 
of orthogonal coordinates such that the coordinate vector fields 
diagonalize the second fundamental form of the immersion at any point. 
This is a strong conclusion since, for instance, it 
allows to express additional assumptions on the submanifold in 
terms of partial differential equations in a system of special 
coordinates.

The attempt in \cite{VZ} to obtain a classification of the Einstein 
submanifolds $f\colon M^4\to\R^6$ with flat normal bundle only 
yielded all possible pointwise structures of the second fundamental 
form. But even in this rather special case there are many 
possibilities for the second fundamental form, as can be seen 
in the Appendix of this paper. Somehow, it is not a surprise 
that there is not a description given, outside of some very simple 
cases, of the actual isometric immersions that carry these second 
fundamental forms.

Warped products of Riemannian manifolds have been intensively 
used to construct interesting classes of Einstein manifolds. 
The class of those that have a surface as base was discussed 
in Chapter $9$ of Besse \cite{Be}. The main purpose of this paper
is to classify, locally or globally, the Euclidean Einstein 
submanifolds in codimension two that belong to this class. 
A classification is also given for higher codimension, but 
now under an extrinsic assumption.

Theorem $9.119$ in \cite{Be}, which is stated without a proof, 
claims a classification of the complete Einstein manifolds that 
are warped products of the form $M^n=L^2\times_\varphi N^{n-2}$, 
$n\geq 4$. It is said that the manifold is either a Riemannian 
product (i.e., the warping function $\varphi$ is constant) 
or it belongs to one of  four families being $L^2=\R^2$ and 
$N^{n-2}$ an Einstein manifold. Some additional information 
about these examples can be seen in the Appendix of \cite{CPW}. 
Although at the present time there is no proof in the 
literature, experts in the field believe that the claim in 
\cite{Be} is true.

\vspace{1ex}
The complete Einstein manifolds that fit our 
purposes are given next.  
\vspace{1ex}

\noindent {\it The  Clifford tori.}
If $\rho>0$ we have the Riemannian product of manifolds
\be\label{Clifford}
M^n=\Sf^2(1/\sqrt{\rho})\times\Sf^{n-2}(\sqrt{(n-3)/\rho})
\ee 
where $\Sf^m(r)$ denotes an $m$-dimensional sphere of 
radius $r$.
\vspace{1ex}

\noindent {\it The Ricci flat Generalized Schwarzschild metric.}
The examples $a)$ in \S $9.118$ of \cite{Be} are of the form
$$
M^n=\R^2\times_\varphi N^{n-2}
$$ 
where $N^{n-2}$ is an Einstein manifold of Ricci curvature 
$\rho=n-3$. The precise description of the Ricci flat metric 
on $M^n$ is given in the next section.
It was called in \cite{Be} and \cite{CPW} the Generalized 
Schwarzschild metric since for $n=4$ and $N^2=\Sf^2(1)$ it 
has been named the Riemannian Schwarzschild metric in analogy 
to the well-known Lorentzian version.

Examples $2$ given in the last section shows that 
$M^n=\R^2\times_\varphi\Sf^{n-2}(1)$ endowed with a 
Generalized Schwarzschild metric can always be isometrically 
immersed in Euclidean space with codimension two as an 
$(n-2)$-rotational submanifold. Recall that an 
\emph{$(n-2)$-rotational submanifold} $f\colon M^n\to\R^{n+p}$, 
$n\geq 3$, with axis $\R^{p+1}$ over a surface 
$g\colon L^2\to\R^{p+2}$ is the $n$-dimensional submanifold 
generated by the orbits of the points of $g(L)$ (disjoint 
from $\R^{p+1}$) under the action of the subgroup $SO(n-1)$ 
of $SO(n+p)$ which keeps pointwise $\R^{p+1}$ invariant.
\vspace{1ex}

Throughout this paper, we denote by $F^m(\varepsilon)$, $m\geq 2$, 
an Einstein manifold with Ricci curvature $(m-1)\varepsilon$, 
$\varepsilon=1,-1,0$, (that is, with normalized Ricci curvature $\varepsilon$). 
This the case of the sphere $\Sf^m(1)$, the hyperbolic space $\Hy^m(-1)$ 
and the Euclidean space $\R^m$.  
\vspace{1ex}

In the realm of complete manifolds the following is the main 
result of the present paper.

\begin{theorem}\label{cod2g}  
Let $M^n=L^2\times_\varphi F^{n-2}(\varepsilon)$, $n\geq 5$, be a complete 
simply connected Einstein manifold whose sectional curvature is 
not constant on any open subset. If $F^{n-2}(\varepsilon)$ 
has constant sectional curvature and $f\colon M^n\to\R^{n+2}$ 
is an isometric immersion, then one of the following holds:
\begin{itemize}
\item[(i)] $M^n=\Sf^2(r_1)\times\Sf^{n-2}(r_2)$ is a Clifford torus 
as in \eqref{Clifford} and $f$ is the product of inclusions 
$\Sf^{n_j}(r_j)\subset\R^{n_j+1}$, $n_j=2,n-2$. 
\item[(ii)]  $M^n=\R^2\times_\varphi\Sf^{n-2}(1)$ is endowed 
with the Generalized Schwarzschild metric and $f$ is an
$(n-2)$-rotational submanifold given by Examples $2$.
\end{itemize}
\end{theorem}

If $M^n$ is not simply connected, then the composition of $f$
with the projection from its universal cover has the same image 
as $f$.
A key ingredient in our proof of the above global theorem is a 
local classification result that characterizes a class of Einstein 
$(n-2)$-rotational submanifolds in codimension two with arbitrary 
Ricci curvature; see our Theorem \ref{cod2} below.
\vspace{1ex}

We are indebted to Andrei Derdzinski, Peter Petersen, 
Thomas Ivey and Wolfgang Ziller for their helpful comments.

\section{A class of Einstein manifolds}

Warped products of Riemannian manifolds have been used to 
construct interesting classes of Einstein manifolds. 
The ones that have a surface as base are discussed in 
Chapter $9$ of Besse \cite{Be}. A few facts from this 
book are exposed in this section.

\begin{proposition}\label{twoeq} 
The warped product $M^n=L^2\times_\varphi F^{n-2}(\varepsilon)$, 
$n\geq 4$, is an Einstein manifold of Ricci curvature $\rho$ 
if and only if the warping function $\varphi\in C^\infty(L)$ 
satisfies the 
equations 
\be\label{eq1}
(n-2)\hess\varphi=(K-\rho)\varphi I
\ee
and 
\be\label{eq2}
\Delta\varphi
+\frac{n-3}{\varphi}(\|\n\varphi\|^2-\varepsilon)+\rho\varphi=0
\ee
where $K$ is the Gaussian curvature of $L^2$ and $\hess\varphi$ 
is the endomorphism of $TL$ determined by the Hessian of $\varphi$. 
Moreover, equations  \eqref{eq1} and \eqref{eq2} for non-constant 
solutions  $\varphi>0$ are equivalent to the equation
\be\label{eq0}
\hess\varphi=\frac{1}{2}(2\Delta\varphi
+\frac{n-3}{\varphi}(\|\nabla\varphi\|^2-\varepsilon)+\rho \varphi)I.
\ee
\end{proposition}

\proof Equations \eqref{eq1}  and \eqref{eq2} are, 
respectively, a particular case of equations $9.107\, b)$ 
and $9.107\,c)$ in \cite{Be}. The equivalence with \eqref{eq0} 
was observed in \S 9.116 of \cite{Be}.
\vspace{2ex}\qed 

The following result is Lemma $9.117$ in \cite{Be} but  
can also be seen in  \cite{Ku} and \cite{Ta}. It states that
the existence of a non-constant solution of equation \eqref{eq0} 
implies that $L^2$ is locally an ``intrinsic surface of rotation", 
namely, it has a warped product structure. 

\begin{proposition}\label{coordinates}
If equation \eqref{eq0} admits a non-constant solution 
$\varphi\in C^\infty(L)$ then there exist local coordinates 
$(t,u)$ on $L^2$ such that $\varphi=\varphi(t)$, the metric of 
$L^2$ has the form $ds^2=dt^2+\varphi'^2 du^2$ and 
$\hess\varphi=\varphi''I$. 
\end{proposition}

Equation \eqref{eq0} in the coordinates given
by Proposition \ref{coordinates} has the form
\be\label{eq}
2\varphi\varphi''+(n-3)(\varphi'^2-\varepsilon)+\rho\varphi^2=0.
\ee
Then, multiplying by $\varphi'\varphi^{n-4}$ and integrating 
gives that
\be\label{firstin}
\varphi'^2=\varepsilon-\frac{\rho}{n-1}\varphi^2+\frac{c}{\varphi^{n-3}}
\ee
where $c\in\R$ is a constant. 
\vspace{1ex}

We summarize the above in the following statement.

\begin{proposition}\label{structure}
Let $M^n=L^2\times_\varphi F^{n-2}(\varepsilon)$ be a warped product 
such that $\n\varphi\neq 0$ at any point. If $L^2$ in the coordinates 
$(t,u)$ carries the metric 
\be\label{metr}
ds^2=dt^2+\varphi'^2 du^2,
\ee 
where $\varphi$ is a solution of \eqref{firstin}, then $M^n$ is an 
Einstein manifold with Ricci curvature $\rho$. Conversely, if $M^n$ 
is an Einstein manifold with Ricci curvature $\rho$ then there 
are local coordinates  $(t,u)$ on $L^2$ such that the metric 
is as in \eqref{metr} where $\varphi$ solves \eqref{firstin}.
\end{proposition}

The Gauss curvature of $L^2$ endowed with the metric 
\eqref{metr} is $K=-\varphi'''/\varphi'$.
It is easily seen that $L^2$ has constant Gauss curvature 
$K=\rho/(n-1)$ if and only if $\varphi$ is a solution of 
equation \eqref{firstin} for $c=0$. If this is the case and 
if $F^{n-2}(\varepsilon)$ has constant sectional curvature then 
it can be easily proved that $M^n=L^2\times_\varphi F^{n-2}(\varepsilon)$ 
has constant sectional curvature $K$. 

\begin{remark}\label{remark}
{\em We observe that equation \eqref{eq} is just equation \eqref{eq2}. 
If willing to verify that any solution of equation \eqref{eq} is also 
a solution of equation \eqref{eq1}, take the derivative of \eqref{eq} 
and then use that $K=-\varphi'''/\varphi'$.
}\end{remark}

We now turn our attention to the case of Einstein manifolds with 
complete metrics in this way wrapping up the discussion started in 
the previous section. Recall that a warped product of two Riemannian 
manifolds is complete if and only if both factors are complete. 
A fact relevant to this paper is that a complete manifold can be 
represented by a warped product where one of the factor is 
not complete. For instance, this is the case of polar 
coordinates for Euclidean space. 
\vspace{1ex}

The warped product $M^n=\R^2\times_\varphi F^{n-2}(1)$ is said 
to be endowed with the Ricci flat Generalized Schwarzschild metric 
when $\R^2$ has the rotationally invariant warped metric 
$ds^2=dt^2+\varphi'^2(t)d\theta^2$ 
in polar coordinates being $\varphi\in C^\infty([0,+\infty))$ 
the unique positive solution of \eqref{firstin} (with $\rho=0$) 
for a given initial condition and a given constant $c<0$.

\section{The local case}

We first show that Einstein manifolds that belong to the class 
discussed above may admit local isometric immersions in Euclidean 
space as $(n-2)$-rotational submanifolds with codimension at least two. 
\vspace{1ex}

\noindent\textbf{Examples 1.} Observe that
an $(n-2)$-rotational submanifold $f\colon M^n\to\R^{n+p}$
can be parametrized as follows.  The manifold $M^n$ is isometric 
to (an open subset of) a warped product 
$L^2\times_\varphi\Sf^{n-2}(1)$ and there is an orthogonal 
splitting $\R^{p+2}=\R^{p+1}\oplus\spa\{e\}$, $\|e\|=1$, such 
that the \emph{profile} $g$ of $f$ has the form $g=(h,\varphi)$
where $h\colon L^2\to\R^{p+1}$ and $\varphi=\<g,e\>>0$. Then  
\be\label{paramet}
f(x,y)=(h(x),\varphi(x)\phi(y))
\ee
where $\phi\colon\Sf^{n-2}(1)\to\R^{n-1}$ denotes the inclusion 
$\Sf^{n-2}(1)\subset\R^{n-1}$.
\vspace{1ex}

\noindent $(\bf{a})$  Given $\rho\in\R$ choose an 
initial condition for $\varphi$ and a value for $c$ such that 
the right hand side of the differential equation \eqref{firstin} 
is positive and less than one. Then,  on some open interval 
$I\subset\R$ there is a solution $\varphi\in C^\infty(I)$ 
of \eqref{firstin} that satisfies $0<\varphi'^2(t)<1$. 
Endow $N^2=I\times I_0$, where $I_0\subset\R$ is also an 
open interval, with the warped metric 
$$
ds^2=(1-\varphi'^2(t))dt^2+\varphi'^2(t)du^2.
$$
Given $p\geq 2$ let $h\colon U\to\R^{p+1}$ be an isometric 
immersion of an open subset $U$ of $N^2$. The metric 
induced by the immersion $g\colon U\to\R^{p+2}$ given by
$$
g(t,u)=(h(t,u),\varphi(t))
$$
is $ds^2=dt^2+\varphi'^2(t)du^2$. 

Now let $f\colon M^n\to\R^{n+p}$ be the immersion 
\eqref{paramet} where $M^n=U\times_\varphi\Sf^{n-2}(1)$.
By Proposition \ref{structure} we have that $f$ is an 
$(n-2)$-rotational Einstein submanifold with Ricci 
curvature $\rho$ whose profile is $g=(h,\varphi)$. 

Now assume that $p=2$. The normal bundles of $h$ and $g$ are 
related by the orthogonal splitting
$$
N_gL=N_hL\oplus\spa\{\delta=(-\varphi' h_t,1-\varphi'^2)\}.
$$
Since
$$
(\nab_{\partial t}\delta)_{T_gL}
=(-\varphi''h_t-\varphi'h_{tt},-2\varphi'\varphi'')_{T_gL}
=-\varphi''g_t
$$
and 
$$
(\nab_{\partial u}\delta)_{T_gL}=(-\varphi' h_{tu},0)_{T_gL}
=-\varphi''g_u,
$$
we have that the second fundamental form of $g$ satisfies 
$A^g_\delta=\varphi''I$, thus $g$ has flat 
normal bundle. On the other hand, from Lemma $2.4$ of \cite{DN} 
(or just by computing the  second fundamental form) we know  
that an $(n-2)$-rotational submanifold has flat normal bundle 
if and only if the profile has this  property.   Hence, in this
case $f$ given by \eqref{paramet} has flat normal bundle.
\vspace{1ex}

\noindent $(\bf{b})$  Assume that $\varphi\in C^\infty(I)$ satisfies 
the stronger condition 
$$
\varphi'^2(t)+\varphi''^2(t)<1,\;\;t\in I.
$$ 
Let $g\colon L^2=I\times[0,2\pi)\to\R^4$ be the 
surface of rotation defined by
$$
g(t,\theta)
=(\psi(t),\varphi'(t)\sin\theta,\varphi'(t)\cos\theta,\varphi(t))
$$
where $\psi\in C^\infty(I)$ is determined by 
$\psi'^2=1-\varphi'^2-\varphi''^2$. Then, we have that 
$f\colon L^2\times_\varphi\Sf^{n-2}(1)\to\R^{n+2}$ given by 
$$
f(t,\theta,y)
=(\psi(t),\varphi'(t)\sin\theta,\varphi'(t)\cos\theta,\varphi(t)\phi(y))
$$
is an $(n-2)$-rotational Einstein submanifold with Ricci 
curvature $\rho$ and flat normal bundle whose profile is $g$. 
\vspace{1ex}

A submanifold $f\colon M^n\to\R^{n+p}$, $n\geq 4$, is 
called \emph{$(n-2)$-umbilical} if its tangent bundle carries 
a maximal $(n-2)$-dimensional umbilical distribution $\U$. 
Thus there exists a smooth non vanishing principal normal 
vector field $\eta$ such that
$$
{\cal U}=\{ X\in TM:\a(X,Y)=\<X,Y\>\eta\;\,
\mbox{for all}\;\,Y\in TM\}
$$
where $\a\colon TM\times TM\to N_fM$ is the second fundamental 
form of $f$.  Since  $\dim{\cal U}\geq 2$ it is well-known that 
$\eta$ is a Dupin principal normal. This means that ${\cal U}$ 
is integrable being the leaves round spheres and 
$\eta$ is parallel in the normal connection along the leaves.
The class of $(n-2)$-umbilical submanifolds has shown up in the 
literature in rather different geometric situations, for instance 
see \cite{AD}, \cite{DF}, \cite{DFT} and \cite{Mo}.
\vspace{1ex}

The $(n-2)$-rotational submanifolds described earlier
are the simplest examples of $(n-2)$-umbilical submanifolds. 
From either Lemma $6$ in \cite{DF} or Proposition $10$ in \cite{DFT} 
we have the following characterization.

\begin{proposition}\label{DF}
Let $f\colon M^n\to\R^{n+p}$ be an $(n-2)$-umbilical 
submanifold. If the orthogonal distribution ${\cal U}^\perp$ 
to the umbilical  distribution ${\cal U}$ is totally geodesic 
in $M^n$ then $f$ is locally an $(n-2)$-rotational submanifold. 
\end{proposition}

Next we characterize a class of Einstein submanifolds with 
arbitrary codimension under an assumption of extrinsic nature.

\begin{theorem}\label{umbilical}
Let $f\colon M^n\to\R^{n+p}$, $n\geq 4$, be an $(n-2)$-umbilical 
isometric immersion of an Einstein manifold of Ricci curvature 
$\rho$ which does not have constant sectional curvature on any 
open subset. Then $p\geq 2$ and each 
point of an open dense subset of $M^n$ has an open neighborhood 
$V\subset M^n$ where one of the following holds:
\begin{itemize}
\item[(i)] $V=U\times W$ is part of a Clifford torus and 
$f|_V=g\times i$ is the product of an isometric immersion 
$g\colon U\subset\Sf^2(1/\sqrt{\rho})\to\R^{p+1}$ and the 
inclusion $i\colon W\subset\Sf^{n-2}(\sqrt{(n-3)/\rho})\to\R^{n-1}$.
\item[(ii)] $f|_V$ is an $(n-2)$-rotational submanifold as 
in part $(a)$ of Examples $1$.
\end{itemize}
\end{theorem}

\proof Since the dimension of the umbilical distribution in 
the assumption is maximal having this property, it follows 
from the Cartan-Fialkow result discussed in the introduction 
that $p\geq 2$. 

We show next that the distribution ${\cal U}^\perp$ is totally 
geodesic in $M^n$. From the Gauss equation, we have
$$
\Ric(X,Y)
=n\<\a(X,Y),H\>-\sum_{i=1}^n\<\a(X,X_i),\a(Y,X_i)\>
$$
where $\{X_1,\ldots,X_n\}$ is an orthonormal frame and
$H=\frac{1}{n}\sum_{i=1}^n\a(X_i,X_i)$ 
the mean curvature vector field of the submanifold.

In the sequel, we denote $\a_{ij}=\a(e_i,e_j)$, $1\leq i,j\leq 2$, 
where $\{e_1,e_2\}$ is an orthonormal frame spanning ${\mathcal{U}^\perp}$.
\vspace{1ex}

\noindent $(i)$ For $X=Y=e_i$, we obtain 
$$
\rho+\|\a_{ii}\|^2+\|\a_{ij}\|^2
=n\<\a_{ii},H\>,\;i\neq j.
$$
Equivalently, we have
$$
\rho+\|\a_{ii}\|^2+\|\a_{ij}\|^2
=\<\a_{ii}, \a_{ii}+\a_{jj}+(n-2)\eta\>,\;i\neq j,
$$
where $\eta$ is the principal normal vector field. Hence
\be\label{ga1}
\rho-K(\mathcal{U}^\perp)=(n-2)\<\a_{ii},\eta\>
\ee
where $K(\mathcal{U}^\perp)$ denotes the sectional curvature.
In particular, we have
\be\label{eqalpha}
\<\a_{11},\eta\>=\<\a_{22},\eta\>.
\ee
Then, since the frame $\{e_1,e_2\}$ is arbitrary, we also have
\be\label{eqalpha2} 
\<\a_{12},\eta\>=0.
\ee

\noindent $(ii)$ For unit vectors $X=Y\in {\cal U}$, we obtain
$$
\rho+\|\eta\|^2=n\<\eta,H\>. 
$$
Equivalently, using \eqref{eqalpha} we have
\be\label{eqalpha1}
\rho-(n-3)\|\eta\|^2=2\<\a_{ii},\eta\>,\;\; 1\leq i\leq 2.
\ee 
Next we use the Codazzi equation
\begin{align*}
\n_X^\perp \a(Y,Z)-\a(\n_X Y,Z)&-\a(Y,\n_X Z)\\
&=\n_Y^\perp \a(X,Z)-\a(\n_Y X,Z)-\a(X,\n_Y Z).
\end{align*}
For $Y=e_i,Z=e_j$ with $ i\neq j$ and $X\in {\cal{U}}$,
we obtain
$$
\n_X^\perp \a_{ij}-\a(\n_X e_i,e_j)-\a(e_i,\n_X e_j)
=\n_{e_i}^\perp \a(X,e_j)-\a(\n_{e_i} X,e_j)-\a(X,\n_{e_i} e_j).
$$
Equivalently, we have
$$
\n_X^\perp \a_{ij}+\<\n_X e_i,e_j\>(\a_{ii}-\a_{jj})
=\< X,\n_{e_i} e_i\>\a_{ij}+\< X,\n_{e_i} e_j\>(\a_{jj}-\eta).
$$
Taking the inner product with $\eta$ and using \eqref{eqalpha}, 
\eqref{eqalpha2} and that $\n_X^\perp \eta=0$ gives
$$
\<\n_{e_i} e_j,X\>(\<\a_{jj},\eta\>-\|\eta\|^2)=0,\;i\neq j.
$$
Then using \eqref{eqalpha1} yields
\be\label{eqeq}
(\rho-(n-1)\|\eta\|^2)\<\n_{e_i}e_j,X\>=0,\;i\neq j.
\ee

We claim that $\rho-(n-1)\|\eta\|^2\neq 0$ on an open dense 
subset. Otherwise, there is an open subset of $M^n$  where 
$\|\eta\|^2=\rho/(n-1)$. We obtain from \eqref{ga1} and 
\eqref{eqalpha1} that
$$
K(\mathcal{U}^\perp)=\frac{\rho}{n-1}\cdot
$$
But this is easily seen to imply that $M^n$ has constant 
sectional curvature on an open subset, and that is a 
contradiction. It follows from \eqref{eqeq} that 
\be\label{codre}
\<\n_{e_i}e_j,X\>=0,\;i\neq j,
\ee
for any $X\in {\cal U}$. 

For $Y=Z=e_i$ and $X\in {\cal U}$, we obtain
$$
\n_X^\perp\a_{ii}-2\a(\n_X e_i,e_i)
=\n_{e_i}^\perp \a(e_i,X)-\a(\n_{e_i} X,e_i)-\a(\n_{e_i}e_i,X).
$$
Equivalently, we have
$$
\n_X^\perp \a_{ii}-2\<\n_X e_i,e_j\>\a_{ij}
=\<\n_{e_i}e_i,X\>(\a_{ii}-\eta)
+\<\n_{e_i}e_j,X\>\a_{ij},\;i\neq j.
$$
Using \eqref{codre} we obtain
$$
\n_X^\perp \a_{ii}-2\<\n_X e_i,e_j\>\a_{ij}
=\<\n_{e_i}e_i,X\>(\a_{ii}-\eta),\;i\neq j.
$$
Taking the inner product with $\eta$ and using \eqref{eqalpha2} 
gives
$$
\<\n_X^\perp \a_{ii},\eta\>
=(\<\a_{ii},\eta\>-\|\eta\|^2)\<\n_{e_i}e_i,X\>,\;\; 1\leq i\leq 2.
$$
On the other hand, from \eqref{eqalpha1} and since 
$\n_X^\perp\eta=0$ we have
$$
\<\n_X^\perp \a_{ii},\eta\>=0,\;\; 1\leq i\leq 2.
$$
Using \eqref{eqalpha1} we obtain
$$
(\rho-(n-1)\|\eta\|^2)\<\n_{e_i}e_i,X\>=0,\;\; 1\leq i\leq 2,
$$
and therefore
\be\label{codre2}
\<\n_{e_i}e_i,X\>=0,\;\; 1\leq i\leq 2,
\ee
for any $X\in {\cal U}$. That the distribution ${\cal U}^\perp$ 
is totally geodesic follows from \eqref{codre} and \eqref{codre2}.

From the above we have that Proposition \ref{DF} applies. Hence
$M^n$ is part of a warped product $L^2\times_\varphi\Sf^{n-2}(1)$.
For simplicity, we assume that either $\varphi$ is constant or 
that $\n\varphi\neq 0$ on $M^n$. In the former case, we obtain 
from \eqref{eq1}  that the Gauss curvature of $L^2$ is the 
constant $K=\rho$.  
Since a Riemannian product of manifolds of constant sectional 
curvature $M_{c_1}^p\times M_{c_2}^{n-p}$ is an Einstein 
manifold if and only if only
\be\label{procon}
(p-1)c_1=(n-p-1)c_2,
\ee
we have from \eqref{paramet} that we are in case $(i)$ of 
the statement. In the latter case, it follows from 
Proposition \ref{structure} that we are in case $(ii)$. \qed

\begin{remark} {\em  If $f\colon M^n\to\R^{n+p}$ is an 
$(n-1)$-umbilical Einstein submanifold then $M^n$ has 
constant sectional curvature. This can be proved by a
simpler version of the above computation or just by using
that a locally conformally flat Einstein manifold has
constant sectional curvature.
}\end{remark}

The following is the main result of this section. Observe 
that only assumptions of intrinsic nature are required.

\begin{theorem}\label{cod2} Let $M^n=L^2\times_\varphi F^{n-2}(\varepsilon)$, 
$n\geq 5$, be an  Einstein manifold such that neither
the warping function nor the sectional curvature are constant
on open subsets of $L^2$ and $M^n$, respectively.
If $F^{n-2}(\varepsilon)$ has constant sectional curvature then 
any isometric immersion $f\colon M^n\to\R^{n+2}$  is an 
$(n-2)$-rotational submanifold as in part $(a)$ of 
\mbox{Examples $1$} along connected components of an open dense 
subset of $M^n$.
\end{theorem}

\begin{remark} {\em Implicit in the above statement there is a
rigidity fact. Namely, any local isometric deformation of the 
submanifold $f$ as above is necessarily the result of a local 
isometric deformation of the surface $h$ in $\R^3$.
}\end{remark}

Let $V^k\subset\R^k$, $k\geq 1$, denote the half space 
$$
V^k=\{x\in\R^k:\sigma(x)=\<x,e\>>0\}
$$
where $e\in\R^k$ is a unit vector. 
A \emph{warped product representation} of the Euclidean 
space $\R^m$, $m=k+\ell$ and $\ell\geq 2$, is the explicitly 
constructible isometry 
$\psi\colon V^k\times_\sigma\Sf^\ell(r)\to\R^m$ onto $\R^m$ 
up to a subspace $\R^{k-1}\subset\R^m$. The warping function 
is a linear function $\sigma(x)=\<x,v\>$ where $\|v\|=1$. By the 
\emph{warped product of isometric immersions} $h_1\colon L^p\to V^k$ 
and $h_2\colon N^q\to\Sf^\ell(r)$ determined by the  warped product 
representation $\psi$ of $\R^m$ we mean the isometric immersion
$f=\psi\circ(h_1\times h_2)\colon L^p\times_{\sigma\circ h_1}N^q\to\R^m$.
Notice that $\sigma\circ h_1\in C^\infty(L^p)$ is just a 
coordinate function of $h_1$.

\begin{lemma}\label{warprod}  
Assume that the warped product $M^n=L^2\times_\varphi N^{n-2}$, $n\geq 5$,
of Riemannian manifolds has the following properties: $(a)$ $L^2$ does 
not contain an open subset of points of constant Gauss curvature and 
$(b)$ $M^n$ does not admit a local isometric immersion into $\R^{n+1}$. 
Given an isometric immersion $f\colon M^n\to\R^{n+2}$, there exists an 
open dense subset of $M^n$ each of whose points lies in an open product 
neighborhood $U=L^2_0\times N^{n-2}_0\subset M^n$ such that $f|_U$ is a 
warped product of isometric immersions with respect to a warped product 
representation $\psi\colon V^{2+k_1}\times_\sigma\Sf^{n-2+k_2}\to\R^{n+2}$, 
$V^{2+k_1}\subset\R^{2+k_1}$, where either $k_1=k_2=1$ or $k_1=2$ and 
$k_2=0$.
\vspace{5ex}

\begin{picture}(150,84)\hspace{-15ex}
\put(100,30){$L^2_0\times_{\sigma\circ h_1} N^{n-2}_0$}
\put(112,55){$h_1$}\put(165,32){\vector(1,0){135}}
\put(305,30){$\R^{n+2}$} \put(154,55){$h_2$}
\put(205,45){$\circlearrowright$}
\put(180,20){$f|_U=\psi\circ(h_1\times h_2)$}
\put(107,42){\vector(0,1){30}} \put(149,42){\vector(0,1){30}}
\put(98,80){$V^{2+k_1}\times_\sigma\Sf^{n-2+k_2}$}
\put(175,80){\vector(3,-1){125}} \put(235,65){$\psi$}
\end{picture}
\end{lemma}

\proof It follows from the assumption $(b)$ that $M^n$ does 
not contain an open subset of flat points. Then Theorem $14$ in 
\cite{DT0} applies and, according to that result, there are 
three disjoint possibilities named there $(i)$ through $(iii)$.
But our assumption $(a)$ excludes the case $k_1=0$ in $(i)$ as well 
as $(iii)$ whereas the assumption $(b)$ excludes $(ii)$.\qed 

\begin{remark}{\em That $k_1=2,k_2=0$  just says that $f|_U$ 
is an $(n-2)$-rotational submanifold with profile $h_1$.
}\end{remark}

\noindent {\em Proof of Theorem \ref{cod2}:}  
We have from Proposition \ref{structure} that $L^2$
admits local coordinates $(t,u)$ such that its metric 
has the form $ds^2=dt^2+\varphi'^2 du^2$ being $\varphi$ a 
solution of equation \eqref{firstin} where $c\neq 0$
since $L^2$ has no constant Gauss curvature. 
Therefore, according to Lemma~\ref{warprod} it remains to show 
that the case $k_1=k_2=1$ there cannot occur. If otherwise, we 
have $h_1=(h^0,\varphi)\colon L^2\to\R^3$ where $h^0=h^0(t,u)$ 
is a local parametrization of $\R^2$ in  an open neighborhood 
of any point in the open dense subset where $\varphi'\neq 0$. 
Then
$$
\|h^0_t\|^2=1-\varphi'^2(t),\;\;\<h^0_t,h^0_u\>=0\;\;\mbox{and}\;\;
\|h^0_u\|^2=\varphi'^2(t).
$$
Thus, there is an orthonormal frame $\{e_1,e_2\}$ such that
$$
h^0_t=\sqrt{1-\varphi'^2}\,e_1\;\;\mbox{and}\;\;h^0_u=\varphi'e_2.
$$
From $\n_{\p u}h^0_t=\n_{\p t}h^0_u$ we derive that 
$e_j=e_j(u),\; j=1,2$, and that
$$
\<\n_{\p u}e_1,e_2\>=\frac{\varphi''}{\sqrt{1-\varphi'^2}}\cdot
$$
It follows that
$$
h^0_{tt}=-\frac{\varphi'\varphi''}{\sqrt{1-\varphi'^2}}e_1
\;\;\;\mbox{and}\;\;\;h^0_{tu}=\varphi''e_2.
$$
Hence
$$
h^0_{ttu}=-\frac{\varphi'\varphi''}{\sqrt{1-\varphi'^2}}\n_{\p u}e_1
=-\frac{\varphi'\varphi''^2}{1-\varphi'^2}e_2
\;\;\;\mbox{and}\;\;\; h^0_{tut}=\varphi'''e_2.
$$
Thus, we have that
$$
\varphi'''+ \frac{\varphi'\varphi''^2}{1-\varphi'^2}=0
$$
or, equivalently, that
$$
(\varphi''/\sqrt{1-\varphi'^2})'=0.
$$
Then $\varphi=(1/a)\sin at$, $0\neq a\in\R$, and this is a 
contradiction since $\varphi$ is not a solution of equation 
\eqref{firstin} with $c\neq 0$.

We have shown that $k_1=2,k_2=0$ and, in particular, that we have  
that $F^{n-2}(\varepsilon)=\Sf^{n-2}(1)$. 
\vspace{2ex}\qed

Theorem \ref{cod2} does not hold if we allow $M^n$ to have 
constant sectional curvature. In this situation, we can still 
exclude the case $(iii)$ in Theorem $14$ of \cite{DT0} 
since the manifold there does not carry a foliation by 
$(n-2)$-dimensional spheres.  Taking a second derivative of 
\eqref{firstin} and using that  $K=-\varphi'''/\varphi'$ it 
follows that $K=\rho/(n-1)$. By Theorem $14$ of \cite{DT0} 
there are two possible situations for which there are 
non-rotational submanifolds. We give next examples for both 
cases in the notation of Lemma \ref{warprod}.
\vspace{1ex}

\noindent $(1)$  We have $L^2=\R^2$ endowed with Cartesian 
coordinates $(t,u)$, $\varphi(t)=t$ and 
$h_2\colon F^{n-2}(1)\to\Sf^n(r)$ 
is any non-totally geodesic submanifold. For instance, 
an appropriate product of spheres. Then $f(M)$ is a Ricci 
flat submanifold but it is not flat.
\vspace{1ex}

\noindent $(2)$ $\R^{n+2}=\R^3\times_\sigma\Sf^{n-1}(r)$ where 
$\sigma=\<x,e\>$ and $r=(n-4)/(n-3)$. 
Let $L^2\subset\R^2$ be endowed with the coordinates $(t,u)$ 
and the metric induced by $h_1\colon L^2\to\R^3$ given below. 
Then let $h_2\colon F^{n-2}(1)\to\Sf^{n-1}(r)$ be a torus being
$$
F^{n-2}(1)=\Sf^m(r_1)\times\Sf^{n-m-2}(r_2),\;\;2\leq m\leq n-4,
$$ 
where $r_1=\sqrt{(m-1)/(n-3)}$ and $r_2=\sqrt{(n-m-3)/(n-3)}$.  
\vspace{1ex}

\noindent Case $K=1$:  Take
$h_1(u,t)=(\cos t\sin u,\cos t\cos u,\sin t)$.
\vspace{1ex}

\noindent Case $K=0$: Take $h_1(u,t)=(\sin u,\cos u,t)$.
\vspace{1ex}

\noindent Finally, we observe that a similar construction for 
$K=-1$ does not work.

\section{The complete case}

We first show that $M^n=\R^2\times_\varphi\Sf^{n-2}(1)$ endowed
with the Generalized Schwarzschild metric can be realized 
parametrically as a rotational submanifold in codimension two  
with flat normal bundle. The case when the dimension is $n=5$ 
turns out to be special in the sense that the parametrization 
is completely explicit.
\vspace{1ex}

\noindent\textbf{Examples 2.} We make use of the construction 
of submanifolds given in part $(b)$ of Examples $1$. 
Let $\varphi\in C^\infty([0,+\infty))$ be the unique positive 
solution of 
$$
\varphi'^2=1+\frac{c}{\varphi^{n-3}}
$$ 
such that  
$\varphi(0)=(n-3)/2$ and $c=-((n-3)/2)^{n-3}$. From
$$
\varphi'^2+\varphi''^2
= 1-((n-3)/2\varphi)^{n-3}+((n-3)/2\varphi)^{2(n-2)}
$$
we have
$$
\varphi'^2(t)+\varphi''^2(t)<1
\;\;\text{if and only if}\;\;\varphi(t)>(n-3)/2.
$$
Since $\varphi'(t)\geq 0$ the inequality holds true for $t>0$. 

We have that $M^n=\R^2\times_\varphi\Sf^{n-2}(1)
=[0,+\infty)\times_{\varphi'}\Sf^{1}(1)\times_{\varphi}\Sf^{n-2}(1)$
is endowed with the Generalized Schwarzschild metric.  Let 
$h\colon\R^2\to\R^3$ be any surface isometric to the surface
of rotation 
$$
h_0(t,\theta)=(\psi(t),\varphi'(t)\sin\theta,\varphi'(t)\cos\theta).
$$
Then the isometric immersion $f\colon M^n\to\R^{n+2}$ given by 
\be\label{expl}
f(t,\theta,y)
=(h(t,\theta),\varphi(t)\phi(y))
\ee
is an $(n-2)$-rotational submanifold with flat normal bundle.
\vspace{1ex}

We recall that by the classical result of Bour any surface of 
rotation in $\R^3$ admits a $2$-parameter family of isometric 
helicoidal surfaces; for instance see Lemma $2.3$ of \cite{CD}.
\vspace{1ex}

For dimension $n=5$ the above yields a totally explicit 
parametrization of the submanifold. In fact, in this case the 
solutions of equation \eqref{eq} are
$$
\varphi(t)=\sqrt{t^2-c},\;\;t^2>c\;\;\mbox{and}\;\;c\in\R. 
$$
In particular, taking $\varphi(t)=\sqrt{t^2+1}$ in \eqref{expl} 
we obtain explicitly parametrized complete Ricci flat Einstein 
submanifolds.

\begin{remark}{\em Solving  \eqref{firstin} for some dimensions, 
namely, for  $n=4,7$ and $9$  leads to elliptic integrals and hence 
the solutions can be expressed in terms  of elliptic functions
whereas in  other cases the integrals are hyperelliptic.
}\end{remark}

\begin{proposition}\label{propglo}
Let  $M^n=L^2\times_\varphi F^{n-2}(1)$, $n\geq 4$, be a 
complete Einstein manifold such that $L^2$ is non-compact and 
$\varphi\in C^\infty(L)$  is non-constant. If there is an 
isometric immersion $g\colon L^2\to\R^m$, $m\geq 4$, such 
that $\varphi=\<g,e\>$ where $e\in\R^m$, $\|e\|=1$, then
the following facts hold:
\begin{itemize}
\item[(i)] $M^n$ is Ricci flat. 
\item[(ii)] If $L^2$ is simply connected then 
$M^n=\R^2\times_\varphi F^{n-2}(1)$ is endowed with the 
Generalized Schwarzschild metric.
\end{itemize}
\end{proposition}

\proof First observe that $\rho\leq 0$. In fact, if otherwise 
$M^n$ is compact, and hence also $L^2$ would be compact which
has been excluded.

Notice that $\n\varphi=e^\top$ gives $\|\nabla\varphi\|\leq 1$.
We argue that $\rho=0$. Taking traces in \eqref{eq1} gives
\be\label{eq3}
(n-2)\Delta \varphi=2(K-\rho)\varphi.
\ee
Then combining with \eqref{eq2} yields
\be\label{eq4}
2K+(n-4)\rho
=(n-2)(n-3)\frac{1-\|\n \varphi\|^2}{\varphi^2}\cdot
\ee
Suppose that $\rho<0$.
We have from \eqref{eq4} that
\be\label{desig}
2K\geq(4-n)\rho.
\ee 
Thus, it follows for $n\geq 5$ that $L^2$ is compact 
and this is a contradiction.

It remains to consider the case $n=4$. We have that 
\be\label{moregen}
\n\|\n\varphi\|^2
=\frac{n-3}{\varphi}(1-\|\n\varphi\|^2)\n\varphi
-\rho\varphi\n\varphi
\ee
on the subset 
$$
L_0=\{x\in L^2:\n\varphi(x)\neq 0\}.
$$
In fact, in the coordinates of Proposition \ref{coordinates}
we have that \eqref{moregen} is equation \eqref{eq}. 
We obtain from \eqref{desig} that $L^2$ is parabolic. It follows 
from \eqref{moregen} that
\be\label{conta}
\nabla(\varphi(\|\n\varphi\|^2-1)+(\rho/3)\varphi^3)=0
\ee
on $L_0$. Therefore, on each connected component of the subset 
$L_0$ we have that
\be\label{neweq}
\|\n\varphi\|^2=1-\frac{\rho}{3}\varphi^2+\frac{C}{\varphi}
\ee
for some constant $C<0$ such that $\varphi^3\leq C/\rho$.
On the other hand, we obtain from \eqref{conta} that the set 
${\rm int}\{x\in L^2:\n\varphi(x)=0\}$ is empty. Hence
\eqref{neweq} holds on all of $L^2$.  In particular, we have 
that $\varphi$ is bounded. It follows from \eqref{eq3} and 
\eqref{desig} that $\Delta\varphi>0$. Since $L^2$ is parabolic 
then $\varphi$ is constant, and this is a contradiction
to \eqref{eq3}.
\vspace{1ex}

Assume that $L^2$ is simply connected. Since $\rho=0$, then 
from \eqref{eq1}, \eqref{eq2}, \eqref{eq3} and \eqref{eq4} 
we have that
\begin{align}
\hess\varphi 
&=\frac{1}{n-2}K\varphi I\label{teq1},\\
\Delta\varphi
&=\frac{n-3}{\varphi}(1-\|\n\varphi\|^2)\geq 0\label{teq2},\\
2K
&= (n-2)(n-3)\frac{1-\|\n\varphi\|^2}{\varphi^2}\geq 0\label{teq4}.
\end{align}
In particular, since $K\geq 0$ then $L^2=\R^2$.
\vspace{1ex}

\noindent\emph{Fact 1:}
There exists a constant $c<0$ such that 
\be\label{teq17}
\|\n\varphi\|^2=1+\frac{c}{\varphi^{n-3}}
\ee
on $L^2$. Moreover, the following equations hold on $L^2$:
$$
(i)\, K= -\frac{(n-2)(n-3)c}{2\varphi^{n-1}},\;
(ii)\,\Delta \varphi=-\frac{(n-3)c}{\varphi^{n-2}},\;
(iii)\,\hess\varphi=-\frac{(n-3)c}{2\varphi^{n-2}}I.
$$

We have from \eqref{moregen} that 
$$
\n(\varphi^{n-3}(\|\n\varphi\|^2-1))=0.
$$
on $L_0$. Therefore, we have that \eqref{teq17} holds for some 
constant on connected components of $L_0$. Moreover, it holds 
trivially on connected components  of 
${\rm int}\{x\in L^2:\n\varphi(x)=0\}$. 
Thus, by continuity \eqref{teq17} holds on all of $L^2$. 

Equations $(i), (ii)$ and $(iii)$ follow from \eqref{teq1}, 
\eqref{teq2}, \eqref{teq4} and \eqref{teq17}. Finally, we 
argue that $c<0$. If otherwise and since $K\geq 0$, we
have from $(i)$ that $c=0$. Then $(i)$ and $(ii)$ give that 
$K=0$ and that $\varphi$ is  a positive harmonic function. 
It follows from Liouville's theorem that 
$\varphi$ is constant, and this is a contradiction due to 
\eqref{teq17}. 
\vspace{1ex}

\noindent\emph{Fact 2:} We have that
$(i)\;\inf \varphi=(-c)^{1/(n-3)}>0,\; (ii)\;\sup\varphi
=+\infty$, $(iii)\; K\; \text{is bounded}$ and $(iv)\;\inf K=0.$
Moreover, a point $p_0$ is a critical point of $\varphi$ 
if and only if 
\be\label{critical}
\varphi(p_0)=\inf \varphi=(-c)^{1/(n-3)}.
\ee

Since $\varphi>0$, by the well-known Omori-Yau maximum principle 
there exists a sequence of points $\{x_m\}_{m\in{\mathbb{N}}}$ 
in $L^2$ such that 
$$
\lim_{m\to \infty} \varphi(x_m)=\inf\varphi,
\;\;\text{and}\;\; \lim_{m\to \infty}\|\n\varphi\|(x_m)=0.
$$
From \eqref{teq17} we have
$$
1+\frac{c}{(\inf\varphi)^{n-3}}=0,
$$
and this is $(i)$. To prove $(ii)$ suppose that $\sup\varphi<+\infty$. 
Since $L^2$ is parabolic and $\varphi$ is subharmonic due to 
part $(ii)$ of Fact $1$, then $\varphi$ is constant, 
a contradiction by $(ii)$ of Fact $1$.  
From $(i)$ of Fact $1$ we have 
$$
0\leq K=-\frac{(n-2)(n-3)c}{2\varphi^{n-1}}
\leq-\frac{(n-2)(n-3)c}{2(\inf\varphi)^{n-1}}
$$
which proves $(iii)$. Since $\sup\varphi=+\infty$ the equality 
part above proves $(iv)$. Finally, it follows from \eqref{teq17} 
that \eqref{critical} holds at a critical point. 
\vspace{1ex}

\noindent\emph{Fact 3:}
The function $\varphi\colon L^2\to (0,+\infty)$ has at most 
one critical point.
\vspace{1ex}

Let $p_0$ be a critical of $\varphi$. By \eqref{critical} 
$p_0$ is a global minimum of $\varphi$. Let $p\neq p_0$
and $\gamma\colon\R\to L^2$ be a geodesic parametrized by 
arc-length such that $\gamma(0)=p_0$ and $\gamma(\ell)=p$. Now
consider the function $u=\varphi\circ\gamma\colon\R\to (0,+\infty)$.
Then $\min_{s\in\R}u(s)=u(0)$. 

From $(iii)$ of Fact $1$ we have that
$$
\n_{\dot{\gamma}}\n u
=-\frac{(n-3)c}{2 u^{n-2}}\dot{\gamma}\cdot
$$
Then
$$
{\dot{\gamma}}\<\n u,\dot{\gamma}\>
=\<\n_{\dot{\gamma}}\n u,\dot{\gamma}\>
=-\frac{(n-3)c}{2u^{n-2}}
$$
and thus
$$
\ddot{u}=-\frac{(n-3)c}{2u^{n-2}}>0.
$$
Thus $\dot{u}$ is monotone increasing which implies for 
$s<0$ that $\dot{u}(s)<\dot{u}(0)=0$ and for $s>0$ that 
$\dot{u}(s)>\dot{u}(0)=0$. Hence $t=0$ is the only point 
where $u$ attains its minimum. Thus $p_0$ is the only point 
of $\gamma$ that $\varphi$ attains its minimum or, equivalently,
the only point of $L^2$ where $\n\varphi$ vanishes. \vspace{2ex}

 Let $p_0\in L^2$ be the necessarily unique possible critical 
point of $\varphi$. Then on either $L^2$ or $L^2\smallsetminus\{p_0\}$ 
we define $e_1=\n \varphi/\|\n \varphi\|$ and $e_2=J(e_1)$ where 
$J$ is the complex structure of $L^2$. From $(iii)$ in Fact $1$, 
we have
$$
X(\|\n\varphi\|)e_1+\|\n\varphi\| \n_X e_1
=-\frac{(n-3)c}{2\varphi^{n-2}}X
$$
for any $X\in TL$. It follows that
$$
\n_{e_1} e_1=0,\;\;
e_1(\|\n\varphi\|)=-\frac{(n-3)c}{2\varphi^{n-2}},\;\;
e_2(\|\n\varphi\|)=0
$$
and
$$
\|\n\varphi\|\<\n_{e_2}e_1,e_2\>=-\frac{(n-3)c}{2\varphi^{n-2}}\cdot
$$
Then
$$
[e_1,\|\n\varphi\|e_2] 
=\|\n\varphi\|(\n_{e_1} e_2-\n_{e_2} e_1)+e_1(\|\n\varphi\|)e_2=0.
$$
Hence, there exist coordinates $(t,u)$ such that 
$$
\d/\d t=e_1\;\;\text{and}\;\; \d/\d u=\|\n\varphi\|e_2
$$
and in these coordinates the metric of $L^2$ is 
$ds^2=dt^2+\|\n\varphi\|^2 du^2$ where
$\frac{\d}{\d u}\|\n\varphi\|^2=0$.
Thus, we can write
\be\label{ODE}
\varphi'^2=1+\frac{c}{\varphi^{n-3}}\cdot
\ee 
Let $\phi_1(x,t)$ and $\phi_2(x,t)$ be the one-parameter groups 
of diffeomorphisms generated by $\n\varphi$ and $J(\n\varphi)$, 
respectively. Since $\|\n\varphi\|=\|J(\n\varphi)\|<1$, it follows 
that $\phi_1(x,t)$ and $\phi_2(x,t)$ are defined for all 
$x\in L^2$ and $t\in\R_+$ or $\R$ according to the existence
or not of a critical point of $\varphi$. Hence $(t,u)$ are either 
Cartesian or polar coordinates in $\R^2$.

To conclude the proof, we argue that the case of Cartesian 
coordinates cannot occur due to the completeness of $L^2$. 
In fact, if these coordinates occur we would have from 
\eqref{ODE} that
$$
\lim_{t\to -\infty}\int_t^{t_0}\frac{\varphi'(u)du}
{\sqrt{1+c\varphi^{3-n}(u)}}=+\infty.
$$
On the other hand, we have
$$
 \int_t^{t_0}\frac{\varphi'(u)du}{\sqrt{1+c\varphi^{3-n}(u)}}
=\int_{\varphi(t)}^{\varphi(t_0)}\frac{dx}{\sqrt{1+cx^{3-n}}}
=A\int_{\sqrt{1+c\varphi^{3-n}(t)}}^{\sqrt{1+c\phi^{3-n}(t_0)}}
\frac{dy}{(1-y^2)^{\frac{n-2}{n-3}}}
$$
where $A=2(-c)^{\frac{1}{n-3}}/(n-3)$ and
$y=\sqrt{1+cx^{3-n}}$, hence $y\in(0,1)$. But then
$$
\lim_{t\to -\infty}\int_t^{t_0}\frac{\varphi'(u)du}
{\sqrt{1+c\varphi^{3-n}(u)}}=
A\int_0^{\sqrt{1+c\phi^{3-n}(t_0)}}
\frac{dy}{(1-y^2)^{\frac{n-2}{n-3}}}<+\infty,
$$
and this is a contradiction.\qed

\begin{corollary}\label{corol}
Let $f\colon M^n\to\R^{n+p}$, $n\geq 4$, be a complete 
$(n-2)$-rotational submanifold. If $M^n$ is an Einstein 
manifold with Ricci curvature $\rho$, then one of the 
following holds:
\begin{itemize}
\item[(i)]  We have that $\rho>0$ and 
$M^n=L^2\times\Sf^{n-2}(\sqrt{(n-3)/\rho})$ where 
$L^2$ is either $\Sf^2(1/\sqrt{\rho})$ or the projective space 
${\R}{\mathbb P}^2(1/\sqrt{\rho})$. 
\item[(ii)] The universal cover $\R^2\times_\varphi\Sf^{n-2}(1)$ 
of $M^n$ is endowed with the Generalized Schwarzschild metric.
\end{itemize}
\end{corollary}

\proof
By assumption we have that $M^n=L^2\times_\varphi\Sf^{n-2}(1)$
where $L^2$ is complete. If $L^2$ is compact then also $M^n$ 
is compact, and from either \cite{Kim} or Theorem $1.2$ of 
\cite{CSW} it follows that $\varphi$ is constant. Then from  
\eqref{eq1} and  \eqref{procon} we obtain $K=\rho>0$, and 
hence $L^2$ is as in part $(i)$ of the statement. 
If $L^2$ is non-compact the proof follows from 
Proposition \ref{propglo}.
\vspace{2ex}\qed

\noindent{\it Proof of Theorem \ref{cod2g}}: 
If $\n\varphi=0$ on an open subset it follows easily
from Theorem 17 in \cite{DT0} that $f$ is locally as in part 
$(i)$, hence it is $(n-2)$-rotational.
On the other hand, from Proposition \ref{coordinates} we have 
that in a neighborhood of any point of $L^2$ where $\n\varphi\neq 0$ 
there are local coordinates $(t,u)$ such that $\varphi=\varphi(t)$ 
satisfies  \eqref{firstin}. It follows from Theorem \ref{cod2} 
that $f$ is locally an $(n-2)$-rotational submanifold along an 
open dense subset $U$ of $M^n$. 
Since the manifold is complete and the dimension of the umbilical 
distribution is minimal, then it is well known (see \cite{Re}) that 
on $U$ the umbilical leaves are complete spheres, and hence $f$ is 
globally an $(n-2)$-rotational submanifold.  The proof now  follows 
from Corollary \ref{corol}.\vspace{2ex}\qed

\begin{example} {\em Let $M^n=\R^2\times_\varphi N^{n-2}$ be
endowed with the Generalized Schwarzschild metric where
$$
N^{n-2}=\Sf^m(r_1)\times\Sf^{n-m-2}(r_2),\;\;2\leq m\leq n-3,
$$ 
for $r_1^2=(m-1)/(n-4)$ and $r_2^2=(n-m-3)/(n-4)$.   
Then  $f\colon M^n\to\R^{n+3}$ given by 
$$
f(t,\theta,y)
=(h(t,\theta),\varphi(t) j(y)),
$$
where $h$ is as in Examples $2$ and $j\colon N^{n-2}\to\R^n$ 
the inclusion, is an isometric immersion that is not rotational.
}\end{example}

\noindent \emph{Problem:} In respect to Theorem \ref{cod2g} notice 
that  for dimension $n=4$ the classification remains open.

\section{Appendix}

Let $f\colon M^4\to\R^6$ be an isometric immersion with flat  normal
of an Einstein manifold.  What makes this case special is that a four 
dimensional Riemannian manifold is Einstein if and only if its sectional 
curvature  tensor satisfies that $K(X,Y)=K(Z,W)$ for any orthonormal 
tangent vectors $X,Y,Z,W$. Besides submanifolds of constant sectional 
curvature or products of space forms, it was observed in \cite{VZ} 
that the possibilities for the second fundamental form of $f$ are as
follows. There exists an orthonormal frame $\{\xi_1,\xi_2\}$ such that 
the corresponding shape operators with respect to a tangent frame that 
diagonalizes the second fundamental form have the form:
$$
A_{\xi_1}=\begin{bmatrix}
a&0&0&0\\
0&b&0&0\\
0&0&c&0\\
0&0&0&d&
\!\!\!\!\end{bmatrix},\;\;\;
A_{\xi_2}=\begin{bmatrix}
0&0&0&0\\
0&p&0&0\\
0&0&q&0\\
0&0&0&r&
\!\!\!\!\end{bmatrix}\\
$$
where $a\neq 0$ and
$$
pq=ad-bc,\;\;\;pr=ac-bd,
\;\;\;qr=ab-cd,
$$
$$
(ba-cd)(ca-bd)(da-bc)>0.
$$
If we drop the last condition, we have
$$
A_{\xi_1}=\begin{bmatrix}
a&0&0&0\\
0&\epsilon a&0&0\\
0&0&b&0\\
0&0&0&\epsilon b&
\!\!\!\!\end{bmatrix},\;\;\;
A_{\xi_2}=\begin{bmatrix}
0&0&0&0\\
0&0&0&0\\
0&0&p&0\\
0&0&0&q&
\!\!\!\!\end{bmatrix}\\
$$
where $a\neq 0$, $pq=\epsilon(a^2-b^2)$ and $\epsilon=\pm 1$.
The submanifolds in the case $\epsilon=1$ are the ones given
by part $(ii)$ of Theorem \ref{umbilical} when $p=2$ and the 
surface $g\colon L^2\to\R^4$ has flat normal bundle. 
On the other hand, making use of the Codazzi equations one
can show that the case $\epsilon=-1$ does not represent the
second fundamental form of any actual submanifold.

\noindent Marcos Dajczer\\
IMPA -- Estrada Dona Castorina, 110\\
22460--320, Rio de Janeiro -- Brazil\\
e-mail: marcos@impa.br
\bigskip

\noindent Christos-Raent Onti\\
Department of Mathematics and Statistics\\
University of Cyprus\\
1678, Nicosia -- Cyprus\\
e-mail: onti.christos-raent@ucy.ac.cy
\bigskip

\noindent Theodoros Vlachos\\
University of Ioannina \\
Department of Mathematics\\
Ioannina -- Greece\\
e-mail: tvlachos@uoi.gr
 
\end{document}